\documentclass[12pt, oneside]{article}

\usepackage{graphicx}
\usepackage{wrapfig}
\usepackage{multicol}
\usepackage{amsthm}
\usepackage{amssymb}
\usepackage{amsmath}
\usepackage{amsfonts}
\usepackage[T1]{fontenc}
\usepackage[english]{babel}
\usepackage{vmargin}

\newtheorem{deft}{Definition}[section]
\newtheorem{pro}{Proposition}[section]
\newtheorem{rem}{Remark}[section]

\newtheorem{coro}{Corollary}[section]
\newtheorem{theo}{Theorem}[section]

\newtheorem{ex}{Example}[section]
\newtheorem{exs}{Examples}[section]
\setmargins{2cm}{1.5cm}{17cm}{24cm}{0cm}{5mm}{1cm}{1.5cm}

\title{Sufficient optimality conditions for a separable product quasiconcave programming 
\thanks{%
Mathematics Subject Classifications: 26B25, 90C26, 46N10.}
}

\author{ M.BERDI \thanks{M.BERDI, Mohammed V University. Rabat. Morocco, medberdi@gmail.com} \, M. A. ABOUCHOUAR \thanks{M. A. ABOUCHOUAR, Mohammed V University. Rabat. Morocco, abou.aminos@gmail.com} \,   A.HASSOUNI \thanks{A.HASSOUNI, Mohammed V University. Rabat. Morocco, hassouni@fsr.ac.ma} }

\date{}
\begin{document}
\parindent=0mm
\maketitle
\bigskip

\begin{small}
\begin{abstract} {\rm }
In mathematical economics, the used functions are, in general, considered to be quasiconcave. Moreover,  they are, in many cases, separable of nature. It is known that a local maximum of a quasiconcave function is not, in general, a global maximum.\\
 In this paper we will show that this property is true when the quasiconcave function is furthermore separable.  Sufficient optimality conditions for a separable quasiconcave programming will be studied in both differentiable and non differentiable cases. Thanks to separability condition, quasiconcave functions have  nice properties in optimization problems.

\end{abstract}
\end{small}

{\textbf{Keywords:} Concave, Quasiconcave, Pseudoconcave Function, Separable  Product Function,  Generalized Superdifferential, modified K.K.T conditions.
\section{Introduction}
Concavity is a central tool in mathematical economics and optimization theory, but in practice the widely used functions in these areas are considered to be quasiconcave instead of concave.
In Arrow-Enthoven $(1961)$ \cite{A-E}, the concave optimization problem was extended to the quasiconcave programming, and sufficient optimality conditions in differentiable case were obtained. Later, several authors have studied the quasiconvex optimality conditions by means of various generalized gradients (see for instance Hassouni \cite{H},  Hiriart-Urruty \cite{HU}, Martinez-Legaz \cite{MLeg}).  

In consummer theory, the functions studied  in many cases are considered to be separable of nature. Sufficient and necessary condition on  separable utility functions to be quasiconcave was given in Yaari \cite{Y}, Debreu and Koopmans \cite{D-K}, Crouzeix and Lindberg \cite{C-L}, Berdi and Hassouni \cite{B-H}.

In this paper, the unconstrained optimization problem
$$ \rm{(USQP)}:\qquad \rm{\max_{x\in X} f(x)}$$
and the constrained optimization problem
$$\rm{(CSQP)}:\qquad  \rm{ \displaystyle{ \max_{x\in X}f(x)  \quad subject\  to\  h_j(x)\geqslant 0 \ (j=1,...,p)  }} $$
are investigated, where 
$\rm{f}$ and $\rm{h_j}\ $ $(j=1,...,p)$ are quasiconcave on $X$, and
$$\rm{f(x) = \prod_{i=1}^{m} f_{i}(x_i), \qquad x=(x_1,...,x_m)\in X=X_1 \times X_2 \times...\times X_m}$$
with $\rm{f_{i}}$ is a positive non constant  real valued function defined on the nonempty open convex set $X_i$ of $\mathbb{R}^{n_i}$.

The scope of the paper is to obtain sufficient optimality conditions for the unconstrained separable quasiconcave problem $ \rm{(USQP)}$, and first order sufficient optimality conditions of K.K.T type for the constrained separable quasiconcave problem  $ \rm{(CSQP)}$ with quasiconcave constraint functions, in both differentiable and non differentiable cases.

The paper is organized as follows: In section $2$, we recall some definitions and properties that will be needed, such as some characterizations of quasiconcavity and pseudoconcavity in the differentiable case.

 In section $3$, we recall the definition and some properties of a multiplicatively separable $\rm{(MS)}$ function. The basic tool is the multiplicative concavity index introduced by Crouzeix and Kebbour in \cite{C-K}. Under the separability condition $\rm{(MS)}$, quasiconcave differentiable functions become more regular (in convex sense), and then nice results of the problems $ \rm{(USQP)}$ and  $ \rm{(CSQP)}$ are obtained.
 
  In section $4$, we define a generalized superdifferential to extend the notion of pseudoconcavity and some of its fundamental properties without differentiability assumption.
  Some useful and important generalizations of pseudoconcavity in non differentiable case were introduced and studied under various assumptions by Aussel \cite{Aus}, Hassouni and Jaddar \cite{H-J}, and Koml\'osi \cite{K} for related works.
  
   Again, under the separability condition  $\rm{(MS)}$, the problems  $ \rm{(USQP)}$ and  $ \rm{(CSQP)}$ will be studied in the non differentiable case. An appropriate variant  of K.K.T conditions of   $ \rm{(CSQP)}$ will be provided.
\section{Preliminaries and notations.}
We recall some definitions and properties that will be needful in the sequel of this paper.
\begin{deft}\label{deft 2.1.}
Let $C$ be a nonempty convex set in $\mathbb{R}^n$ and let $g$ be a real valued function defined on $C$.\\

The hypograph  of $g$ is the set
$$hyp(g):=\left\lbrace (x,\mu)\in C\times \mathbb{R}:g(x)\geqslant \mu \right\rbrace$$

For $\alpha\in \mathbb{R}$, the upper-level set $U(g,\alpha)$  and the strict upper-level set $U^{s}(g,\alpha)$ of $g$ are defined as follows:
$$U(g,\alpha):=\left\lbrace x\in C: g(x)\geqslant \alpha \right\rbrace.$$
$$U^{s}(g,\alpha):=\left\lbrace x\in C: g(x) > \alpha \right\rbrace.$$
     \begin{itemize}
 \item   $g$ is  said to be \textbf{ concave} on $C$ if $hyp(g)$ is convex, or equivalently if
  $$g((1-\lambda)x+\lambda y) \geqslant (1-\lambda)g(x)+\lambda g(y) $$
 for all $x,y\in C$ and for all $\lambda \in \left[0;1\right].$

 \item $g$ is  said to be \textbf{ quasiconcave} on $C$ if the upper-level set $U(g,\alpha)$ is convex for all $\alpha \in \mathbb{R}$, or equivalently if
  $$g((1-\lambda)x+\lambda y) \geqslant \min\left[g(x);g(y)\right] $$
 for all $x,y\in C$ and for all $\lambda \in \left[0;1\right].$

\item $g$ is  said to be \textbf{ semi-strictly quasiconcave} on $C$ if  
 $$ \forall x_1, x_2 \in C, x_1 \not = x_2 :\ \  g(x_2) > g(x_1) \Rightarrow g((1-\lambda) x_ 1+  \lambda x_ 2) > g(x_ 1) $$
  for all $ \lambda $ in $ (0,1) $

  For a full description of concavity and quasiconvexity we refer to \cite{C-M,C-F,D-A-Z,R}.
\item   $g$ is said to be convex (quasiconvex, semi-strictly quasiconvex) if $-g$ is concave (quasiconcave, semi-strictly quasiconcave).  
  
\item  If $g$ is positive, it is said to be \textbf{logarithmically concave} (log-concave for short) if
   $\ln \circ g$ is  concave.
   \end{itemize}
   The following  properties are rather direct consequences of the definitions.
   
\begin{itemize}
   \item $g$ concave $\Longrightarrow g$ log-concave  $\Longrightarrow g$ quasiconcave.\\
   Let   $\varphi:g(C)\rightarrow \mathbb{R}$;  
  \item If $g$ is quasiconcave on $C$ and  $\varphi$
    is nondecreasing then   $\varphi \circ g$ is  quasiconcave. 
  \item If $g$ is convex (concave) and $\varphi$ is nondecreasing and convex (concave) then $\varphi \circ g$ is convex(concave).
   \item If $g$ is convex (concave) and $\varphi$ is nonincreasing and concave (convex) then $\varphi \circ g$ is concave (convex).\\
 In particular, if $g$ is positive and concave on $C$, then $\displaystyle \frac{1}{g}$ is convex on $C$.    
\item  Given  $x \in C$ and $d\in \mathbb{R}^{n}$, let us   define 
    $$I_{x,d} := \{t\in \mathbb{R} : x+td \in   C \}$$
    $$ g_{x,d}(t)=g(x+td), \qquad t\in I_{x,d} $$
  Then, $g$ is concave (quasiconcave, log-concave) if and only if for every 
  $x\in C$ and $d\in \mathbb{R}^n$,  the function  $g_{x,d}$ is  concave (quasiconcave, log-concave) on the interval $I_{x,d}.$
 \end{itemize} 
\end{deft}
Let's recall the well known characterization of  quasiconcavity under differentiability assumption.
%
%
\begin{pro}\label{pro 2.1.}
(Arrow-Enthoven \cite{A-E})
 Let $g$ be a differentiable function defined on a nonempty open convex set $X$ of $\mathbb{R}^n$.  
Then, $g$ is quasiconcave on $X$ if and only if for all $x_1,x_2 \in X:$
$$\displaystyle{\left\langle \nabla g(x_1),x_2-x_1 \right\rangle< 0 \Rightarrow g(x_2)< g(x_1) }$$
\end{pro}
Under differentiability assumption, the following proposition gives a necessary condition of quasiconcavity by means of the strict upper level-set.
%
%
\begin{pro}\label{pro 2.2.}
Let $g$ be a differentiable quasiconcave function defined on a convex set $X$ of $\mathbb{R}^n$, and $x,\bar{x}\in X$. Then,
$$\displaystyle{x\in clU^s(g, g(\bar{x})) \Rightarrow \left\langle \nabla g(\bar{x}), x-\bar{x}\right\rangle \geq 0 }$$
\end{pro}
By extending the  inequality in Proposition \ref{pro 2.1.}, the notion of  pseudoconcavity  was introduced as a generalized concavity which plays an important role in applied mathematics  such as, optimization theory and mathematical economics. Let's recall the definition  and some properties of pseudoconcavity.
%
\begin{deft}\label{deft 2.2.}
Let $g$ be a differentiable function defined on an open nonempty convex set $X$ of $\mathbb{R}^n$.  \\
$g$ is said to be pseudoconcave on $X$ if for all $x_1,x_2 \in X$
$$\displaystyle{ \left\langle\nabla g(x_1), x_2-x_1 \right\rangle \leq 0  \Rightarrow g(x_2) \leq g(x_1)   }$$
$g$ is said to be pseudoconvex on $X$ if $(-g)$ is pseudoconcave.
\end{deft}
%
\begin{pro}\label{pro 2.3.}(Mangasarian \cite{M})
Let $g$ be differentiable on $X$. Then, 
if $ g $ is pseudoconcave on $ X $, then it is semi-strictly quasiconcave on $X$.
\end{pro}
%
\begin{pro}\label{pro 2.4.}
Let $g$ be a differentiable function on an open convex set $X$ of $\mathbb{R}^n$. Then, $g$ is pseudoconcave on $X$ if and only if the restriction of $g$ to any line segment in $X$ is pseudoconcave. (see \cite{D-A-Z}).
\end{pro}
It has been known that a differentiable pseudoconcave function $g$ is quasiconcave and has a maximum at $x$ whenever $\nabla g(x) =0$. In \cite{C-F}, Crouzeix and Ferland have shown that this property is  a necessary and sufficient condition for pseudoconcavity.
%
%
\begin{pro}\label{pro 2.5.}(Theorem 2.2. \cite{C-F})
Let $g$ be a differentiable and quasiconcave function on an open convex set $X \subset \mathbb{R}^n$ . Then $g$ is pseudoconcave on $X$ if and only if $g$ has a local maximum at $x\in X$ whenever
 $\nabla g(x) = 0$.
\end{pro}
%
\begin{deft}
Let $f$ be a real function defined on a convex set $C \subset \mathbb{R}^n$. For $r\in \mathbb{R}\setminus \left\lbrace0\right\rbrace$, let $f_r$ be defined by: \qquad $f_r(x)=e^{rf(x)}$.\\
$f$ is said to be $r$-concave if $f_r$ is concave whenever $r>0$ and convex whenever $r<0$.\\
$f$ is $0$-concave if it is concave. (For more details on $r$-concavity/convexity see Avriel \cite{Avr}).
\end{deft}
\begin{pro}\label{pro 2.6.}(Theorem 6.1. \cite{Avr})
Let $r$ be any real number and let $f$ be a differentiable 
$r$-concave function on a convex set $C \subset \mathbb{R}^n$. Then $f$ is pseudoconcave on $C$.
\end{pro}

\section{Optimality conditions under Separability: The differentiable case.}
In this section and the next one,  we will study  sufficient optimality conditions for separable quasiconcave  programming when the objective function is multiplicatively decomposed.\\
First, we recall the definition and some properties of a separable product function.
%
%
\begin{deft}\label{deft 3.1.}
Let $X$ be a subset of $\mathbb{R}^n$. 
A function $f:X \rightarrow \mathbb{R}$ is said to be  multiplicatively  separable if it satisfies the following condition:\\
$(MS):$ there exist subsets $X_i$ of $\mathbb{R}^{n_i}$, 
 and functions $f_i:X_i \rightarrow \mathbb{R}$ $\left( i=1,...,m \right)$ such that 
            $$  f(x)=\prod_{i=1}^m f_i(x_i) $$
 where $\displaystyle{ x=(x_1,...,x_m)\in X=X_1 \times \cdots\times X_m}$ 
  and  $\displaystyle{\sum_{i=1}^m n_i=n}$  \\
\end{deft}

First, we recall a necessary condition for the function $f$ to be quasiconcave. See \cite{B-H}
%
%
\begin{pro}\label{pro 3.1.}(Lemma 3.4. \cite{B-H})
Let $f$ be a real function defined on a convex set $X \subset \mathbb{R}^n$ verifying the condition $(MS)$.
If $f$ is quasiconcave on $X$, then all $f_i \ (i=1,...,m)$ are quasiconcave on $X_i$. 
\end{pro}
%
%
Now we recall the definition and some properties of the multiplicative concavity index of a function introduced by Crouzeix and Kebbour in \cite{C-K}. Such  an index was the basic tool to study the (generalized)concavity of a function.
\begin{deft}   {(Crouzeix-Kebbour \cite{C-K})}\label{deft 3.2.}
The multiplicative concavity index $i_{cv}(f)$ of  a function $f:X\rightarrow (0,\infty)$   is defined as follows:
     \begin{itemize}
           \item If there exists $\lambda < 0$ such that $f^{\lambda}$ is not convex then
              $$i_{cv}(f) = \displaystyle \sup \{\mu < 0 : f^{\mu} \ \textrm{ is  convex} \}$$
            \item If $f^{\lambda} $\ is convex for every $\lambda < 0 $, then
              $$i_{cv}(f) = \displaystyle\sup\{\mu > 0 : f^{\mu}\  \textrm{is  concave} \}$$
     \end{itemize}
\end{deft}
The following proposition is an immediate consequence of Definition \ref{deft 3.2.}. (See  \cite{C-K,C-L}).
%
%
\begin{pro}\label{pro 3.2.}Let $X$ and $f$ as in Proposition \ref{pro 3.1.}. Then,
       \item[(a)] If   $\displaystyle {i_{cv}(f) > -\infty}$, then $f$ is quasiconcave  ;   
       \item[(b)] $f$ is log-concave if and only if \  $\displaystyle{ i_{cv}(f) \geq 0}$ ;   
       \item[(c)] $f$ is concave if and only if \ $\displaystyle{i_{cv}(f) \geq 1}$ ;   
       \item[(d)] $f$ is constant if and only if \  $\displaystyle{i_{cv}(f) = +\infty}$ and $X$ is open ;   
       \item[(e)] Let $\displaystyle{\alpha > 0}$, then \  $\displaystyle{i_{cv}(f^{\alpha}) = \frac{1}{\alpha} i_{cv}(f)}$ ;   
       \item[(f)]  $\displaystyle{i_{cv}(f) = \inf \{i_{cv}(f_{x,d}) : x\in X ,d \in \mathbb{R}^{n} \setminus \{ 0 \} \}}$.
\end{pro}

The following proposition,  which reveals an interesting property of separable quasiconcave functions, will be used frequently in the sequel of this paper (see \cite{B-H}).
%
%
\begin{pro}\label{pro 3.3.}(Theorem 3.12.\cite{B-H}) 
For $i = 1,...,m$, let $X_{i}$ be a non-empty open convex subset of $\mathbb{R}^{n_{i}}$, and let $f_{i}$ be a positive non constant function defined on $X_{i}$, and  $f$ be the function defined on the product space $X = X_{1} \times X_{2} \times ... \times X_{m}$ by $$f(x_{1} ,...,x_{m}) =  \displaystyle\prod_{i=1}^{m}f_{i}(x_{i})$$
   \begin{itemize}
      \item[i)] The function $f$ is quasiconcave if and only if one of the following holds:
               \begin{itemize}
                   \item[a)] all functions $f_{i}$ are log-concave.
                   \item[b)] all functions $f_{i}$ except one are log-concave and
                            \begin{equation}\label{eq 1}
                                     \displaystyle\sum_{i = 1}^{n} \frac{1}{i_{cv}(f_{i})} \leq 0
                            \end{equation}
               \end{itemize}
       \item[ii)] If $f$ is quasiconcave then
                        \begin{equation}\label{eq 2}
                                   \frac{1}{i_{cv}(f)} = \displaystyle\sum_{i = 1}^{n} \frac{1}{i_{cv}(f_{i})}
                        \end{equation}
   \end{itemize}
   with the convention $\displaystyle{ \frac{1}{0} = \infty }$. 
\end{pro}
%
%
\begin{rem}\label{rem 3.1}
Notice that  if all $f_{i} \ (i=1,\cdots,m) $ are differentiable log-concave, then so is the function $f$, and then it is pseudoconcave. If not, there exists $i_0\in \left\lbrace1,\cdots,m\right\rbrace$ such that $f_{i_o}$  is not log-concave. It is clear that  $f_{i_0}$ is $i_{cv}(f_{i_0})$-concave with  $i_{cv}(f_{i_0})<0$, then, from Proposition \ref{pro 2.6.}.,   $f_{i_0}$ is pseudoconcave on $X$, hence,    $f$ is also pseudoconcave on $X$.(For more details see the proof of Theorem $5.3.$ in \rm{Crouzeix-Hassouni} \cite{C-H} ).
\end{rem}

\subsubsection*{Unconstrained problem.}
%
%
Consider the unconstrained  problem
$$\displaystyle{(USQP)}:\qquad \displaystyle{ \max_{x\in X} f(x)} $$
where $f$ satisfies the condition $(MS)$ in Definition \ref{deft 3.1.}.
%
%
\begin{theo}\label{Theo 3.1.}
 Assume that  all $f_i$ are differentiable, and $f$ is quasiconcave on $X$. Then,
\item[i)] If $\bar{x}$ is a critical point of $f$ then $\bar{x}$ is a global maximum of $(USQP)$. 
 \item[ii)] If $\bar{x}$ is a local maximum of $(USQP)$ then it is a global maximum.
\end{theo}
%
%
\begin{proof}
 $i)$ Suppose that $\bar{x}$ is a critical point, that is $\nabla f(\bar{x})  =0$, since $\displaystyle{ \left\langle\nabla f(\bar{x}), x-\bar{x}\right\rangle=0 \ \forall x\in X}$ then, by the  pseudoconcavity  of $f$ , one has $ f(x)\leq f(\bar{x})  $ for all $x\in X.$ \\
$ii)$  Let $ \bar{x} $ be a local maximum of $f$, then there is a neighbourhood $\mathcal{N}(\bar{x}) $   of $  \bar{x} $, such that for all $ x \in  \mathcal{N} ( \bar{x})  \cap X $ we have $ f(x) \leq f (\bar{x})$. \\
Let $ x  \in X $ such that $ x  \not  \in  \mathcal{N} ( \bar{x})   $. There exists $  \tilde{ \lambda}  \in (0,1) $ such that $  \tilde{x} = (1-  \tilde{ \lambda})  \bar{x} +  \tilde{ \lambda} x  \in  \mathcal{N} ( \bar{x})  \cap X $. \\
By Proposition \ref{pro 2.4.}. the restriction of $f$ to the line segment $\left[\bar{x},x \right]$ is pseudoconcave. 
If $ f (x) > f( \bar{x}) $ then by the semi-strict quasiconcavity of $ f $ (see Proposition \ref{pro 2.3.}.),  one has $ f ( \tilde {x}) > f( \bar{x }) $ which is a contradiction. Thus $  \bar{x} $ is a global maximum .
\end{proof}
%
%
\subsubsection*{Constrained problem.}
Now consider the constrained  problem
\begin{center}
$(CSQP)$:
\begin{tabular}{lcr}
                 &   $\displaystyle{  \max_{x\in X} f(x)}$    &            \\
subject to       &   $h_j(x)\geq 0$                           &   $j=1,\cdots, p$
\end{tabular}
\end{center}
 where $f$ satisfies the condition $(MS)$ in Definition \ref{deft 3.1.}.\\
Define the feasible set $F=\left\lbrace x\in X :h_j(x) \geq 0, \ j=1,...,p \right\rbrace $.
%
%
\begin{theo}\label{Theo 3.2.}
Assume that $f$ and all  $h_j \ (j=1,...,p)$ are differentiable and quasiconcave on $X$.
Let $\bar{x}$ be a feasible point such that $\displaystyle{ \nabla h_j(\bar{x}) \neq 0  }$ for all $j$ and $h_j(\bar{x}) >0$ for some $j$.\\ If there exist $\lambda_j\in \mathbb{R}, j = 1,\cdots,p,$  such that
\begin{equation}\label{eq 3}
\nabla f (\bar{x})+\sum_{j=1}^p\lambda_j \nabla h_j(\bar{x})=0, 
\end{equation}
\begin{equation}\label{eq 4}
\lambda_jh_j(\bar{x})=0,\  j=1,\cdots,p,  
\end{equation}
\begin{equation}\label{eq 5}
\lambda_j\geq 0,\  j=1,\cdots,p, 
\end{equation}
then $\bar{x}$ is a global solution  of $(CSQP).$
\end{theo}
%
%
\begin{proof}
Assume, by contradiction, that there exists  $x_0 \in F$ such that
$f(x_0) > f(\bar{x})$. By the pseudoconcavity of $f$ we have $\left\langle\nabla f(\bar{x}), x_0- \bar{x}\right\rangle > 0$.\\
Since $h_j(x_0)\geq 0 = h_j(\bar{x}),\ j\in J(\bar{x})=\left\lbrace j / h_j(\bar{x})=0\right\rbrace$, the quasiconcavity of $h_j$ implies $\left\langle\nabla h_j(\bar{x}), x_0-\bar{x}\right\rangle\geq 0,\ j \in J(\bar{x})$.\\
From the complementarity condition (\ref{eq 4}), we have
$\lambda_j = 0, \forall j \not\in J(\bar{x})$, hence:\\ $\displaystyle{\left\langle\nabla f(\bar{x}),  x_0-\bar{x}\right\rangle +\sum_{j=1}^p\left\langle\lambda_jh_j(\bar{x}),  (x_0-\bar{x}\right\rangle > 0}$, 
which contradicts (\ref{eq 3}).
\end{proof}
%
%
\begin{ex}
Consider the Cobb-Douglas utility function defined by:
$$  u(x_{1},x_{2},...,x_{n}) = \prod_{i=1}^{n}u_{i}(x_{i})  $$
where $u_{i}(x_{i}) = x_{i}^{\alpha_{i}} $ with $x_{i} > 0$ and $\alpha_{i} > 0$ ; $i=1,...,n$ \\
Consider the problem of maximization of $u$ with constraint budget:
\begin{center}
$(P)$:
\begin{tabular}{lcr}
            &   $\displaystyle{\max u(x)}$    &    \\
subject to  &   $G(x) \leq B$                 &
\end{tabular}
\end{center}
where $G(x)= p_1x_1+p_2x_2+...+p_nx_n $.

From Proposition \ref{pro 3.3.}.,  $u$ is  quasiconcave , and since it is differentiable with  $\nabla u(x) \neq~0$ for all $x$, then,  by Proposition \ref{pro 2.5.}.,  $u$ is pseudoconcave.\\
Let $H(x)=B-p_1x_1-p_2x_2-...-p_nx_n$. It is clear that 
$H$ is quasiconcave and differentiable.\\
 Let $\bar{x}=(\bar{x}_1,...,\bar{x}_n)$ a feasible point, that is $H(\bar{x})\geq0$, and let $\lambda \in \mathbb{R}$ such that $(\bar{x},\lambda)$ satisfies the $KKT$ conditions (\ref{eq 3}), (\ref{eq 4}) and (\ref{eq 5}) in Theorem \ref{Theo 3.2.}. \\
By (\ref{eq 3}) one has:
$ \displaystyle{ \alpha_i \bar{x}_1^{\alpha_1}\bar{x}_2^{\alpha_2}\cdots \bar{x}_i^{\alpha_i -1}\cdots \bar{x}_n^{\alpha_n}=\lambda p_i }$ for all $i=1 \cdots n$.

Since $\bar{x}_i>0, \alpha_i>0, p_i>0, (i=1,...,n)$ then $\lambda \neq 0$, and $\displaystyle{\frac{\displaystyle{\prod_{j=1}^{n}}(\bar{x}_{j})^{\alpha_{j}}}{\lambda}=\frac{p_i}{\alpha_{i}}\bar{x}_i}$, 
 thus $\displaystyle{  \bar{x}_i= \frac{\alpha_i p_1}{\alpha_1 p_i}\bar{x}_1 }$ for all $i=1,...,n$.\\
By (\ref{eq 4}) and since $\lambda \neq 0$ one has $H(\bar{x}_1,...,\bar{x}_n)=0$, \rm{i.e.} $p_1\bar{x}_1+...+p_n\bar{x}_n=B$,
 then;\\ $\displaystyle{ \bar{x}_1=\frac{\alpha_1 B}{(\alpha_1+...+\alpha_n)p_1} }$,
  thus $\displaystyle{ \bar{x}=\big(\frac{\alpha_1 B}{(\alpha_1+...+\alpha_n)p_1},...,\frac{\alpha_n B}{(\alpha_1+...+\alpha_n)p_n} \big)}$  is a solution of $(P)$.
\end{ex}


\section{Optimality conditions under Separability: The non differentiable case.}

In this section we will study the problems $\rm{(USQP)}$ and $\rm{(CSQP)}$ studied in the previous section when the objective function and the constrained functions are not necessarily differentiable.

First, notice that in Proposition \ref{pro 3.3.}., any assumption of differentiability of $f$ is required. \\
Secondary, we recall that  the notion of pseudoconcavity can be extended to non differentiable case by means of  a generalized superdifferential instead of the classical gradient.\\
Let's  define an abstract superdifferential  in the same sense as the abstract subdifferential defined by Aussel et al in \cite{A-C-L}.
%
%
\begin{deft}\label{deft 4.1}
We call superdifferential, denoted by $\stackrel{\frown}{\partial}$, any operator which associates
a subset $\stackrel{\frown}{\partial}f(x) $ of $\mathbb{R}^n$ to any upper semi-continuous function $f : X \rightarrow \mathbb{R}\cup \left\lbrace -\infty\right\rbrace$
and any $x\in X$ such that $f(x)$ is finite, and satisfies the following properties:
\item[(P1)]$\stackrel{\frown}{\partial} f(x) = \left\lbrace x^*\in \mathbb{R}^n:\left\langle x^*,y-x\right\rangle + f(x) \geq f(y), \forall y \in X \right\rbrace $ whenever $f$ is concave;
\item[(P2)]$0 \in \stackrel{\frown}{\partial}f(x) $ whenever $x$ is a local maximum of $f$;
\item[(P3)]$ \stackrel{\frown}{\partial}(f+g)(x)\subset \stackrel{\frown}{\partial} f(x) + \stackrel{\frown}{\partial} g(x)$ whenever $g$ is a real-valued concave continuous function which is $\stackrel{\frown}{\partial}$-differentiable at $x$.\\
where $g$ is $\stackrel{\frown}{\partial}$-differentiable at $x$ means that both $ \stackrel{\frown}{\partial}g(x)$ and $\stackrel{\frown}{\partial}(-g)(x)$ are
non-empty. 
\end{deft}
%
%
\begin{exs}
Let's recall the Clarcke-Rockafellar subdifferential $\partial^{CR}$ and the upper-Dini subdifferential $\partial^{D+}$ for a lower-semicontinuous function
   $f : X \rightarrow \mathbb{R}\cup \left\lbrace +\infty\right\rbrace$:
 $$\displaystyle{\partial^{CR}f(x):=\left\lbrace x^*\in X^* : \left\langle x^*,v\right\rangle\leq f^{\uparrow} (x;v)  \ \forall v\in X\right\rbrace}   $$
  $${\rm with}\qquad \displaystyle{ f^{\uparrow}(x;v) = \sup_{\varepsilon>0}\limsup_{x'\rightarrow_f x\atop {t\searrow 0}}\inf_{v'\in B_{\varepsilon}(v)}\frac{f(x'+tv')-f(x')}{t},}$$
   and
  $$\displaystyle{\partial^{D^+}f(x):=\left\lbrace x^*\in X^* : \left\langle x^*,v\right\rangle\leq f^{D^+}(x,v),\ \forall v\in X\right\rbrace}$$ 
 $${\rm with}\qquad\displaystyle{f^{D^+}(x,v)=\limsup_{t\searrow 0}\frac{f(x+tv)-f(x)}{t}}$$
 where the notation $\displaystyle{x'\rightarrow_f x} $ means that $x'\rightarrow x$ and $f(x')\rightarrow f(x)$.\\
Let's define the  Clarcke-Rockafellar superdifferential $\stackrel{\frown}{\partial}^{CR}$ and the upper-Dini superdifferential $\stackrel{\frown}{\partial}^{D+}$ for an upper-semicontinuous function  $f : X \rightarrow \mathbb{R}\cup \left\lbrace -\infty\right\rbrace$:  
$$\stackrel{\frown}{\partial}^{CR}f(x):=-\partial^{CR}(-f)(x) $$
$$\stackrel{\frown}{\partial}^{D+}f(x):=-\partial^{D+}(-f)(x)$$
These two superdifferentials check, among others, the above abstract superdifferential's properties  in Definition \ref{deft 4.1}.
\end{exs}
We recall that $\partial^{CR}$ and $\partial^{D+}$  contain the best known subdifferentials such as the lower Hadamard subdifferential $\partial ^{H-}$, the Fr\'echet subdifferential $\partial ^F$ and the lipschitz subdifferential $\partial ^{LS}$.\\
In the sequel, we will use the symbol $\stackrel{\frown}{\partial}$ to mean either $\stackrel{\frown}{\partial}^{CR}$ or $\stackrel{\frown}{\partial}^{D+}$.\\
The following proposition, which extend Proposition \ref{pro 2.1.}., is an immediate consequence of  Definition \ref{deft 4.1}.   and Theorem $2.1.$ in \cite{Aus}.
%
%
\begin{pro}\label{pro 4.1}
 Let $X$ be a nonempty convex subset of $\mathbb{R}^n $ and let $f: X \longrightarrow \mathbb{R}$ be an  upper-semicontinuous  function. Then, the following assertions are
equivalent:
 \item[i)] $f$ is quasiconcave;
 \item[ii)] $\displaystyle{ \left(\exists x^*\in  \stackrel{\frown}{\partial}f(x) : \left\langle x^*,y-x\right\rangle < 0\right) \Rightarrow f(y)  > f(z) \qquad \forall z\in \left[x;y\right])} $
\end{pro}
From Definition \ref{deft 4.1}. and Proposition \ref{pro 4.1}., the  Proposition \ref{pro 2.2.}. can be extended as follows:
%
%
\begin{pro}\label{pro 4.2.}
Let $f$ be an upper semi-continuous and  quasiconcave real valued function on a convex set $X$ of $\mathbb{R}^n$, and let $x,\bar{x}\in X$. Then,
$$x\in clU^s(f,f(\bar{x})) \Rightarrow \left\langle \bar{x}^*,x-\bar{x}\right\rangle \geq 0 \qquad \forall \bar{x}^*\in \stackrel{\frown}{\partial}f(\bar{x})$$ 
\end{pro}
%
Now, we extend the definition and some properties of pseudoconcavity to the non differentiable case.
%
%
\begin{deft}\label{deft 4.2}  Let $X$ be a nonempty convex subset of $\mathbb{R}^n $ and let $f: X \longrightarrow \mathbb{R}$ be an  upper-semicontinuous  function. 
 $f$ is said to be pseudoconcave  with respect to $ \stackrel{\frown}{\partial}$ (in short $\stackrel{\frown}{\partial}$-pseudoconcave) if, for any $x, y\in X$, one has
$$\displaystyle{ \left(  \exists x^*\in  \stackrel{\frown}{\partial} f(x) : \left\langle x^*,y-x\right\rangle \leq 0\right) \Rightarrow f(y)\leq f(x)}.$$
\end{deft}

As in the differentiable case, every $\stackrel{\frown}{\partial}$-pseudoconcave  function  satisfies the following fundamental properties:
%
%
\begin{pro}\label{pro 4.3.}
Let $f$ be an upper-semicontinuous function on $X$. Then, $f$ is $\stackrel{\frown}{\partial}$-pseudoconcave  on $X$  if and only if the restriction of $f$ to any line segment in $X$ is $\stackrel{\frown}{\partial}$-pseudoconcave.
\end{pro}

The proof follows from the radial property  of the generalized derivatives $f^{\uparrow}$ and $f^{D+}$.
%
%
\begin{pro}\label{pro 4.4.}
Let $f: X\rightarrow\mathbb{R}$ be an upper-semicontinuous and radially continuous function. Then, the following assertions are equivalent:
\item[(i)]$f$ is $\stackrel{\frown}{\partial}$-pseudoconcave;
\item[(ii)]$f$ is quasiconcave and ($0\in  \stackrel{\frown}{\partial} f(x) \Rightarrow f$ has a global maximum at $x$).
\end{pro}

The proof is a direct consequence of Definition \ref{deft 4.2}.  and  Proposition \ref{pro 4.1}.                                     
%
%
\subsubsection*{Unconstrained problem}
Let $f$ be a real valued  function defined on a non empty open convex set  $X$ of $\mathbb{R}^n$ satisfying the separability condition $(MS)$ in Definition \ref{deft 3.1.}. Consider the unconstrained problem:
$$\displaystyle{ \rm{(USQP)}: \qquad \max_{x\in X} f(x) }$$
%
%
\begin{theo}\label{theo 4.1.}
 Assume that $f$ is upper-semicontinuous, quasiconcave and radially continuous on $X$. Then,
 if $\bar{x}$ solves $(USQP)$ locally, then  it is a global solution.
 \end{theo}
\begin{proof}
By the separability condition $(MS)$ and from Proposition \ref{pro 3.3.}. and Remark \ref{rem 3.1}., the function  $f$ is either log-concave or $r$-concave, and then it is  $\stackrel{\frown}{\partial}$-pseudoconcave. If $\bar{x}$ is a local maximum of $f$, then by $(P2)$ in Definition \ref{deft 4.1}., $\displaystyle{0\in  \stackrel{\frown}{\partial}f(\bar{x})}$, thus from Proposition \ref{pro 4.4.}., $\bar{x}$ is a global maximum of $f$. 
\end{proof}
%
%
\begin{coro}
Let $f$ as in Theorem \ref{theo 4.1.}. Assume that $\displaystyle{\bar{x}\notin clU^s(f,f(\bar{x}))}$. Then $\bar{x}$ is a global solution of $(USQP)$.
\end{coro}
%
%
\begin{proof}
If $\displaystyle{\bar{x}\notin clU^s(f,f(\bar{x}))}$, then there exists a neighbourhood $\displaystyle{\mathcal{N}(\bar{x})}\subset X$ of $\bar{x}$ such that $\displaystyle{ \mathcal{N}(\bar{x})\cap  clU^s(f,f(\bar{x}))=\emptyset}$, thus $\bar{x}$ is a local maximum of $f$, hence by the $\stackrel{\frown}{\partial}$-pseudoconcavity of $f$ and Theorem \ref{theo 4.1.}., $\bar{x}$ is a global solution of $(USQP)$.
\end{proof}

\subsubsection*{Constrained problem.}
Let $f$ be an upper-semicontinuous and quasiconcave real valued  function defined on $X$ satisfying the separability condition $(MS)$ in Definition \ref{deft 3.1.}. Consider the constrained problem:
\begin{center}
$(CSQP)$:
\begin{tabular}{lcr}
                 &   $\displaystyle{  \max_{x\in X} f(x)}$    &            \\
subject to       &   $h_j(x)\geq 0$                           &   $j=1,\cdots, p$
\end{tabular}
\end{center}
where $h_j \ (j=1,...,p)$ are quasiconcave and upper-semicontinuous on $X$.\\
Let's define the feasible set $\displaystyle{\mathcal{F}=\left\lbrace x\in X :h_j(x) \geqslant 0, \ j=1,...,p \right\rbrace }$.\\
 For $x\in \mathcal{F}$, denote $\displaystyle{ J(x)=\left\lbrace j:h_j(x)=0 \right\rbrace  }$.
%
%
\begin{pro}\label{pro 4.5.}
Let $\bar{x}$ be a feasible point of $(CSQP)$. For  $j\in J(\bar{x})$, assume that $\displaystyle{0 \notin \stackrel{\frown}{\partial}h_j(x)}$ whenever  $x\in \mathcal{F}$ and $h_j(x)=0$. Then,  $\displaystyle{ \left\langle \bar{x}^*,x-\bar{x}\right\rangle \geq 0 \ \forall \bar{x}^*\in \stackrel{\frown}{\partial}h_j(\bar{x})}$.
\end{pro}
%
%
\begin{proof}
For  $j\in J(\bar{x})$, and  $x\in \mathcal{F}$ such that $h_j(x)=0$, if  $\displaystyle{0 \notin \stackrel{\frown}{\partial}h_j(x)}$, then, by $(P_2)$ in Definition \ref{deft 4.1}., $x$ is not a local maximum of $h_j$. Suppose that there exists $\displaystyle{\bar{x}^*\in \stackrel{\frown}{\partial}h_j(\bar{x})} $, such that  $\displaystyle{ \left\langle \bar{x}^*,x-\bar{x}\right\rangle < 0 }$. Since $\displaystyle{h_j(x)=h_j(\bar{x})=0}$,  Proposition \ref{pro 4.2.}. yields  $\displaystyle{x\notin clU^s(h_j, h_j(x))}$. Thus $x$ is a local maximum of $h_j$, which is a contradiction.
\end{proof}

Let's define a variant of the well known K.K.T. conditions that we show as a sufficient optimality  conditions for the constrained problem $(CSQP)$. 
%
%
\begin{deft}\label{deft 4.3.}
We say that a pair $(\bar{x},\bar{\lambda}) \in X \times \mathbb{R}^p$ satisfies the
modified Karush-Kuhn-Tucker conditions (m-K.K.T. conditions) if it satisfies  the super-gradient condition:
\begin{equation}\label{eq 6}
0\in \stackrel{\frown}{\partial}f(\bar{x})+\sum_{j=1}^p\lambda_j \stackrel{\frown}{\partial} h_j(\bar{x})-N_X(\bar{x}) 
\end{equation}
where $\displaystyle{N_X(\bar{x})}$ is the normal cone of $X$ at $\bar{x}$, 
and also  the usual complementary slackness conditions:
\begin{equation}\label{eq 7}
\lambda_jh_j(\bar{x})=0,\  j=1,\cdots,p
\end{equation}
\begin{equation}\label{eq 8}
h_j(\bar{x})\geq 0,\   j=1,\cdots,p
\end{equation}
\begin{equation}\label{eq 9}
\lambda_j\geq 0,\  j=1,\cdots,p
\end{equation}

\end{deft}
Let's recall the {\itshape{Slater constraint qualification}}: there exists  $\tilde{x}$ in $X$, called a Slater point for $\rm{(CSQP)}$, such that 
  $h_j(\tilde{x})>0$ for some $\displaystyle{j\in \left\lbrace1,\cdots,p\right\rbrace}$. 
%
%
\begin{theo}\label{Theo 4.2.}
Let $\bar{x}$ be a feasible point of $(CSQP)$. Assume that $(CSQP)$ has a Slater point, and $0 \not\in \stackrel{\frown}{\partial}h_j(x)$ whenever $x\in \mathcal{F}$ and $h_j(x)=0 $. If there exists $\bar{\lambda}\in \mathbb{R}^p$   such that $(\bar{x},\bar{\lambda})$ satisfies the m-K.K.T. conditions,  then $\bar{x}$ is a solution of $(CSQP)$.
\end{theo}
%
%
\begin{proof}
Assume, by contradiction, that there exists a feasible point $x_0$ such that $f(x_0) > f(\bar{x})$.\\
 By Proposition \ref{pro 3.3.}, because of separability condition, $f$ is actually $\stackrel{\frown}{\partial}$-pseudoconcave  and then, for all $\displaystyle{\bar{x}^{*}\in\stackrel{\frown}{\partial}f(\bar{x})}$ one has  $\displaystyle{\left\langle \bar{x}^*, x_0- \bar{x}\right\rangle > 0}$.\\
 Since $N_X(\bar{x})$ coincides with the normal cone of convex analysis when $X$ is convex (see \cite{Cl}), then for all $v\in N_X(\bar{x})$ one has 
$\left\langle v,x_0-\bar{x}\right\rangle \leq 0$. Thus, for all $\bar{x}^*\in \stackrel{\frown}{\partial} f (\bar{x})$, $\bar{x}_j^*\in \stackrel{\frown}{\partial} h_j(\bar{x})$ and $v\in N_X(\bar{x})$ one has:
\begin{equation}\label{eq 10}
\left\langle \bar{x}^*-v , x_0-\bar{x}   \right\rangle > 0
 \end{equation}
 If $j\in J(\bar{x})$, by Proposition \ref{pro 4.5.}., one has for all $\bar{x}_j^* \in \stackrel{\frown}{\partial} h_j(\bar{x})$:   
  \begin{equation}\label{eq 11}
  \left\langle \bar{x}_j^*, x_0-\bar{x}\right\rangle \geq 0
\end{equation}   
From the  condition (\ref{eq 7}),   
$\lambda_j = 0$ for all $j \not\in J(\bar{x})$.
  Adding (\ref{eq 11}) for $j=1,\cdots,p$, and combining with (\ref{eq 10}) we get: 
   $$\displaystyle{ \left\langle \bar{x}^*, x_0- \bar{x}\right\rangle + \sum_{j=1}^p \lambda_j\left\langle \bar{x}_j^*, x_0-\bar{x}\right\rangle - \left\langle v,x_0-\bar{x} \right\rangle >0}$$
    for all  $\displaystyle{\bar{x}^*\in  \stackrel{\frown}{\partial} f (\bar{x})}$, $\displaystyle{  \bar{x}_j^* \in \stackrel{\frown}{\partial} h_j(\bar{x}), \ (j=1,...,p) }$ and $\displaystyle{v\in N_X(\bar{x})}$, which contradicts (\ref{eq 6}).
\end{proof}



\begin{thebibliography}{99}
\bibitem{A-E}
 K. J. Arrow and A. C. Enthoven, "Quasi-concave programming", Econometrica., 29(1961),
779-800.

\bibitem{Aus}
D. Aussel, Subdifferential properties of quasiconvex and pseudoconvex functions: Unified approach, J. Optim. Theory Appl., 97(1998), 29--45.

\bibitem{A-C-L}
D. Aussel, J. N. Corvellec and M. Lassonde, Mean value property and subdifferential
criteria for lower semicontinuous functions, Trans. Amer. Math. Soc., 347(1995), 4147--4161.

\bibitem{Avr} M.Avriel, r-Convex Functions, Math. Program. 2(1972), 309-323.

\bibitem {B-H}
 M. Berdi and A. Hassouni, Quasiconcavity of a separable product of utility functions, Optimization, 67(11)(2018), 1837-1848.

\bibitem{C-M}
 A. Cambini and L. Martein, Generalized convexity and optimization, Lecture Notes in Econom. and Math. Systems., Springer-Verlag Berlin,  2009.
 
 \bibitem{Cl}
F. H. Clarke, Optimization and nonsmooth analysis, Wiley-Interscience, New York, 1983.

\bibitem{C-F}
 J. P. Crouzeix and J. A. Ferland, Criteria for quasiconvexity and pseudoconvexity :relationships and comparisons, Math. Program.,  23(1982), 193--205.

\bibitem{C-H}
 J. P. Crouzeix and A. Hassouni, Generalized monotonicity of a separable product of operators: the multivalued case, Set-Valued Anal., 3(1995), 351--373.

\bibitem{C-K}
 J. P. Crouzeix and R. Kebbour, Index multiplicatifs de convexit\'e/concavit\'e et applications, Cahiers Centre \'Etudes Rech. Op\'er., 34(1)(1992), 7--23.

\bibitem{C-L}
 J. P. Crouzeix and P. O. Lindberg, Additively decomposed quasiconvex functions, Math. Program., 1986;   35(1986), 42--57.

\bibitem{D-K}
 G. Debreu and T. C. Koopmans, Additively decomposed quasiconvex functions, Math. Program.,   24(1982), 1--38.

 \bibitem{D-A-Z}
 W. E. Diewert, M. Avriel and I. Zang, Nine kinds of quasiconcavity and concavity, J. Econom. Theory.,  25(1981), 397--420.

\bibitem{H}
A. Hassouni,  Sous-Differentiels des Fonctions Quasiconvexes, PhD Thesis,
Universit\'e Paul Sabatier, Toulouse, France, (1983).

\bibitem{H-J}
A. Hassouni and A. Jaddar, On pseudoconvex functions and applications to global optimization, ESAIM: Proceedings, vol.20(2007), 138--148.

\bibitem{HU}
J-B. Hiriart-Urruty, On optimality conditions in nondifferentiable programming, Math. Program. 14(1978), 73--86.

\bibitem{K}
S. Koml\'osi, Some properties of nondifferentiable pseudoconvex functions, Math. Program. 26(1983), 232--237.  

\bibitem{M}
O. L. Mangasarian,  Pseudoconvex-functions, Journal of Siam Control Series A, 3 (1965),  281--290.

\bibitem{MLeg}
J-E. Martinez-Legaz, Weak lower subdifferentials and applications, Optimization, 21(1990), 321--341.

\bibitem{R}
 R. T. Rockafellar, Convex analysis, Princeton University Press, Princeton, New Jersey, 1970.

\bibitem{Y}
 M. E. Yaari, A note on separability and quasiconcavity,  Econometrica., 45(5)(1977), 1183--1186.

\end{thebibliography}
\end{document}